\documentclass[leqno,12pt]{article}  
\usepackage{amsmath}
\usepackage{amssymb}

\usepackage[german,english]{babel}
\usepackage{amsthm}
\usepackage{xypic}  

\author{\textsc{Elmar Grosse-Kl\"onne}}
\date{}

\theoremstyle{plain} 
\newtheorem{satz}{Theorem}[section]  
\newtheorem{kor}[satz]{Corollary}  
\newtheorem{lem}[satz]{Lemma}  
\newtheorem{pro}[satz]{Proposition}  
 
\theoremstyle{remark}

\theoremstyle{definition}

\newcommand{\0}{\ensuremath{\overrightarrow{0}}}

\begin{document}

%

\begin{center}{\bf $p$-torsion coefficient systems for ${\rm SL}_2({\mathbb
      Q}_p)$ and ${\rm GL}_2({\mathbb
      Q}_p)$}\\Elmar Grosse-Kl\"onne
\end{center}
\begin{abstract} We show that the categories of smooth ${\rm SL}_2({\mathbb Q}_p)$-representations (resp. ${\rm GL}_2({\mathbb Q}_p)$-representations) of level $1$ on $p$-torsion modules are equivalent with certain explicitly described equivariant coefficient systems on the Bruhat-Tits tree; the coefficient system assigned to a representation $V$ assigns to an edge $\tau$ the invariants in $V$ under the pro-$p$-Iwahori subgroup corresponding to $\tau$. The proof relies on computations of the group cohomology of a compact open subgroup group $N_0$ of the unipotent radical of a Borel subgroup.  
\end{abstract}

\section{Introduction}

In the seminal paper \cite{ss}, P. Schneider and U. Stuhler explained how the smooth representation theory over ${\mathbb C}$ of a $p$-adic reductive group $G$ can be embedded into the theory of $G$-equivariant coefficient systems on the Bruhat-Tits building $X$ of $G$. A smooth $G$-representation $V$ on a ${\mathbb C}$-vector space is 'spread out' onto $X$, by assigning to a simplex $\tau$ of $X$ the space of invariants of $V$ under a suitable compact open subgroup of $G$ fixing $\tau$. Under suitable hypotheses, the $0$-th homology group of the coefficient system so obtained gives back $V$. One obtains equivalences between categories of $G$-representations generated by their invariants under sufficiently small compact open subgroups of $G$, and certain categories of $G$-equivariant coefficient systems on $X$. The interest in this is that the investigation of coefficient systems, more 'local' objects than $G$-representations, can often be reduced to the representation theory of finite groups. It is known that these principles work effectively also for smooth $G$-representations on vector spaces over fields of positive characteristic different from $p$.

On the other hand, it is quite nebulous if there is a similarly tight relationship between smooth $G$-representations on vector spaces over fields of characteristic $p$ on the one hand, and coefficient systems on $X$ on the other hand. For general $G$ the functors analogous to those in \cite{ss} do not establish analogous equivalences of categories. The purpose of this note is to show that, nevertheless, for the groups $G={\rm SL}_2({\mathbb Q}_p)$ and ${\rm GL}_2({\mathbb Q}_p)$ the theory, when restricted to smooth representations of level $1$, develops parallel to that in \cite{ss}, although with quite different proofs. 

Let $\Lambda$ be an Artinian local algebra with residue class field of characteristic $p$. Let ${\mathfrak R}(\Lambda)$ (resp. ${\mathfrak R}_{\ast}(\Lambda)$) denote the category of smooth $G$-representations (resp. ${\rm GL}_2({\mathbb Q}_p)$-representations) on $\Lambda$-modules which are generated by their invariants under a pro-$p$-Iwahori subgroup. We will describe a concrete and explicit category of $G$-equivariant (resp. ${\rm GL}_2({\mathbb Q}_p)$-equivariant) coefficient systems ${\mathfrak C}(\Lambda)$ (resp. ${\mathfrak C}_{\ast}(\Lambda)$) on the Bruhat-Tits tree $X$ of $G={\rm SL}_2({\mathbb Q}_p)$ (which is the same as the Bruhat-Tits tree of ${\rm PGL}_2({\mathbb Q}_p)$). For both $?=$(the empty symbol) and $?=\ast$ we will also describe a functor ${\mathfrak R}_?(\Lambda)\to{\mathfrak C}_?(\Lambda), V\mapsto {\mathcal F}_V$ such that for an edge $\tau$ of $X$ the value ${\mathcal F}_V(\tau)$ is the submodule of invariants in $V$ under the pro-$p$-Iwahori subgroup $U_{\tau}$ corresponding to $\tau$. Our main result is the following:

\begin{satz}\label{einlmain} (= Corollary \ref{main}) The categories ${\mathfrak R}_?(\Lambda)$ and ${\mathfrak C}_?(\Lambda)$ are abelian. The functors $${\mathfrak R}_?(\Lambda)\longrightarrow{\mathfrak C}_?(\Lambda),\quad\quad V\mapsto {\mathcal F}_V$$ and $${\mathfrak C}_?(\Lambda)\longrightarrow{\mathfrak R}_?(\Lambda),\quad\quad {\mathcal F}\mapsto H_0(X,{\mathcal F})$$ are exact and quasi-inverse to each other.
\end{satz}

We expressively point out the following two consequences.

(i) We have $H_1(X,{\mathcal F})=0$ for each ${\mathcal F}\in {\mathfrak C}_?(\Lambda)$. Therefore, Theorem \ref{einlmain} tells us in particular that for each $V\in {\mathfrak R}_?(\Lambda)$ the chain complex (augmented by $V$)\begin{gather}0\longrightarrow C_1(X,{\mathcal F}_V)\longrightarrow C_0(X,{\mathcal F}_V)\longrightarrow V\longrightarrow0\label{stapre}\end{gather}is exact, i.e. is a $G$-equivariant (resp. ${\rm GL}_2({\mathbb Q}_p)$-equivariant) resolution of $V$. If $V$ is admissible then the $C_j(X,{\mathcal F}_V)$ are finitely generated. The existence of 'standard presentations' (Colmez), i.e. of finitely generated two-term resolutions of admissible ${\rm GL}_2({\mathbb Q}_p)$-representations $V$ has been known before, and it plays an important role in Colmez' recent proof of the $p$-adic local Langlands correspondence for ${\rm GL}_2({\mathbb Q}_p)$. Besides encompassing also the group $G={\rm SL}_2({\mathbb Q}_p)$, our result seems to be new in that it provides {\it canonical} and even {\it functorial} such resolutions.

(ii) Theorem \ref{einlmain} tells us in particular that for each ${\mathcal F}\in {\mathfrak C}_?(\Lambda)$ and each edge $\tau$ of $X$ the natural map ${\mathcal F}(\tau)\to H_0(X,{\mathcal F})^{U_{\tau}}$ is bijective. In other words, we have explicit control of the ${U_{\tau}}$-invariants in $H_0(X,{\mathcal F})$; in general this is the main difficulty in trying to understand the $0$-th homology group assigned to an equivariant coefficient system.

The proof of Theorem \ref{einlmain} is elementary, e.g. compared to known proofs for the existence of 'standard presentations'. First we reduce to the case where $\Lambda$ is itself a field $k$ of characteristic $p$. Then the proof is based on computations of the cohomology, with values in certain $k[N_0]$-modules, for a compact open subgroup group $N_0$ of the unipotent radical of a Borel subgroup $N$ of $G$. As $G={\rm SL}_2({\mathbb Q}_p)$ we have $N_0\cong{\mathbb Z}_p$.

It would be very interesting to know wether Theorem \ref{einlmain} can be extended to higher levels, i.e. to representations generated only by their invariants under higher congruence subgroups.\\

{\it Notations:} Let $p$ be a prime number and let $G={\rm SL}_{2}({\mathbb Q}_p)$. Let $X$ be the Bruhat-Tits tree of $G$. For $j\in\{0,1\}$ denote by $X^j$ its set of $j$-simplices. We generally identify a $j$-simplex with its set of 
vertices. For $x\in X^0$ let $K_x\subset G$ denote the maximal compact
subgroup fixing $x$. For $\tau\in X^1$ let $U_{\tau}$ be the maximal
pro-$p$-subgroup of $G$ fixing $\tau$; this is a pro-$p$-Iwahori subgroup.

Let $NT$ be a Borel subgroup of $G$ with unipotent radical $N$ and
split maximal torus $T$. We fix an edge $\sigma=\{x_+,x_-\}\in X^1$ contained in the
apartment corresponding to $T$. We may assume that the vertices $x_+$, $x_-$ are labelled
in such a way that $N_0:=K_{x_+}\cap N$ fixes $x_-$. We let $X_{+}$ denote the
maximal connected full closed subcomplex of $X$ with $x_{+}\in X_{+}^0$ but
$x_{-}\notin X_{+}^0$. Here $X_{+}^j$, for $j\in\{0,1\}$, denotes the set of $j$-simplices of
$X_+$. The group $N_0$ acts on the half tree $X_{+}$. 

Let ${{\mathfrak o}}$ be a complete discrete valuation ring with residue class field $k$ of characteristic $p$. Let ${\rm Art}({\mathfrak o})$ denote the category of Artinian local ${\mathfrak o}$-algebras with residue class field $k$.\\

{\bf Definition:} Let $\Lambda\in{\rm Art}({\mathfrak o})$. A homological coefficient
system ${\mathcal F}$ in $\Lambda$-modules on $X$ is a collection of data as follows:

--- a $\Lambda$-module ${\mathcal F}(\tau)$ for each simplex $\tau$ 

--- a $\Lambda$-linear map $r^{\tau}_x:{\mathcal F}(\tau)\to {\mathcal F}(x)$ for each $x\in X^0$ and $\tau\in X^1$ with $x\in\tau$.

We obtain a $\Lambda$-linear map \begin{gather}\bigoplus_{\tau\in X^1}{\mathcal F}(\tau)\longrightarrow \bigoplus_{x\in X^0}{\mathcal F}(x)\label{difdef}\end{gather} sending $y\in{\mathcal F}(\tau)$ to $\sum_{x\in X^0\atop x\in\tau} r^{\tau}_x(y)$. The cokernel of the map (\ref{difdef}) is denoted by $H_0(X,{\mathcal F})$, its kernel is denoted by $H_1(X,{\mathcal F})$.\\

We say that the homological coefficient system ${\mathcal F}$ is $G$-equivariant if in addition we are given a $\Lambda$-linear map $g_{\tau}:{\mathcal F}(\tau)\to {\mathcal F}(g\tau)$ for each $g\in G$ and each simplex $\tau$, subject to the following conditions: 

(a) $g_{h\tau}\circ h_{\tau}=(gh)_{\tau}$ for each $g,h\in G$ and each simplex $\tau$ 

(b) $1_{\tau}={\rm id}_{{\mathcal F}(\tau)}$ for each simplex $\tau$ 

(c) $r^{g\tau}_{gx}\circ g_{x}=g_{\tau}\circ r_{x}^{\tau}$ for each $g\in G$ and each $x\in X^0$ and $\tau\in X^1$ with $x\in\tau$.\\

It is clear that if ${\mathcal F}$ is a $G$-equivariant homological coefficient system, then $G$ acts compatibly on the source and on the target of the map (\ref{difdef}), hence it acts on $H_0(X,{\mathcal F})$ and on $H_1(X,{\mathcal F})$. There is an obvious notion of a morphism ${\mathcal F}\to{\mathcal G}$ between $G$-equivariant homological coefficient systems: a collection of maps ${\mathcal F}(\tau)\to{\mathcal G}(\tau)$ for all simplices, compatible with the restriction maps and the $G$-actions.\\
 
Similarly we define ${\rm GL}_{2}({\mathbb Q}_p)$-equivariant homological coefficient systems on $X$.\\

{\bf Definition:} (i) For $\Lambda\in{\rm Art}({\mathfrak o})$ let ${\mathfrak
  C}(\Lambda)$ denote the category of $G$-equivariant homological coefficient
systems ${\mathcal F}$ in $\Lambda$-modules on $X$ satisfying the following
three conditions:\\(a) for any $\tau\in X^1$ the action of $U_{\tau}$ on
${\mathcal F}(\tau)$ is trivial,\\(b) for any $z\in\tau\in X^1$ the map
${\mathcal F}(\tau)\to{\mathcal F}(z)$ is injective and its image is
${\mathcal F}(z)^{U_{\tau}}$, and\\(c) for any $z\in X^0$ the $K_z$-representation ${\mathcal F}(z)$ is generated by ${\mathcal F}(z)^{U_{\tau}}$ for any $\tau\in X^1$ with $z\in\tau$. Explicitly: if $S=\{\tau\in X^1\,|\,z\in\tau\}$, then ${\mathcal F}(z)=\sum_{\tau\in S}{\rm im}[{\mathcal F}(\tau)\to{\mathcal F}(z)]$.

(ii) Let ${\mathfrak R}(\Lambda)$ denote the category of smooth $G$-representations on $\Lambda$-modules which are generated by their ${U_{\sigma}}$-invariants. \\

For $V\in{\mathfrak R}(\Lambda)$ define the coefficient
system ${\mathcal F}_V$ on $X$ by putting ${\mathcal
  F}_V(\tau)=V^{U_{\tau}}$ for $\tau\in X^1$, and ${\mathcal
  F}_V(x)=\sum_{\tau}V^{U_{\tau}}$ (sum inside $V$) for $x\in X^0$, where $\tau$ runs
through the edges containing $x$. The transition map $r^{\tau}_x:{\mathcal F}_V(\tau)\to {\mathcal F}_V(x)$ for $x\in X^0$ and $\tau\in X^1$ with $x\in\tau$ is defined as follows: if $x$ lies in the $G$-orbit $Gx_+$ of $x_+$ then $r^{\tau}_x$ is the inclusion; if however $x\notin Gx_+$ then $r^{\tau}_x$ is {\it the negative of} the inclusion (note that for each $\tau\in X^1$ exactly one vertex of $\tau$ belongs to $Gx_+$). It is clear that ${\mathcal
  F}_V\in{\mathfrak C}(\Lambda)$. Conversely, for ${\mathcal F}\in{\mathfrak
  C}(\Lambda)$ we have $H_0(X,{\mathcal F})\in{\mathfrak R}(\Lambda)$.

\begin{lem}\label{adjoin} These functors form an adjoint pair: for ${\mathcal F}\in {\mathfrak C}(\Lambda)$ and $V\in{\mathfrak R}(\Lambda)$ we have a natural isomorphism$${\rm Hom}_{{\mathfrak R}(\Lambda)}(H_0(X,{\mathcal F}),V)\cong{\rm Hom}_{{\mathfrak C}(\Lambda)}({\mathcal F},{\mathcal F}_V).$$
\end{lem}

{\sc Proof:} First assume that we are given an element $\rho\in{\rm Hom}_{{\mathfrak R}(\Lambda)}(H_0(X,{\mathcal F}),V)$. Given $\tau\in X^1$, let $x\in X^0$ denote the unique vertex with $x\in\tau$ and contained in $Gx_+$. The composition $\rho'_{\tau}:{\mathcal F}(\tau)\to{\mathcal F}(x)\to H_0(X,{\mathcal F})\stackrel{\rho}{\to} V$ is $U_{\tau}$-equivariant, hence its image lies in $V^{U_\tau}={\mathcal F}_V(\tau)$ as $U_{\tau}$ acts trivially on ${\mathcal F}(\tau)$. We obtain a map $\rho_{\tau}:{\mathcal F}(\tau)\to{\mathcal F}_V(\tau)$. It can alternatively be described as follows. Let $y\in X^0$ denote the unique vertex with $x\in\tau$ and which is {\it not} contained in $Gx_+$. Then $\rho_{\tau}'$ is the {\it negative} of the composition ${\mathcal F}(\tau)\to{\mathcal F}(y)\to H_0(X,{\mathcal F})\stackrel{\rho}{\to} V$: this follows from the fact that $\rho$ is well defined modulo the image of the map (\ref{difdef}). To define $\rho_x:{\mathcal F}(x)\to {\mathcal F}_V(x)$ for vertices $x$ we use that ${\mathcal F}(x)$ is the sum of all ${\mathcal F}(\tau)$ with $x\in\tau$. Namely, given $x\in X^0$ contained in $Gx_+$, we define $\rho_x$ as the sum of the maps $\rho'_{\tau}$ with $\tau$ running through all $\tau\in X^1$ with $x\in\tau$. On the other hand, given $x\in X^0$ not contained in $Gx_+$, we define $\rho_x$ as the sum of {\it the negatives of} the maps $\rho'_{\tau}$ with $\tau$ running through all $\tau\in X^1$ with $x\in\tau$; by our second characterization of the maps $\rho'_{\tau}$ given above this is well defined. Together we obtain the element in ${\rm Hom}_{{\mathfrak C}(\Lambda)}({\mathcal F},{\mathcal F}_V)$ which we want to assign to $\rho$. Conversely, let $\xi\in{\rm Hom}_{{\mathfrak C}(\Lambda)}({\mathcal F},{\mathcal F}_V)$ be given. Summing the maps ${\mathcal F}(x)\to{\mathcal F}_V(x)\subset V$ for all $x\in X^0$ provided by $\xi$ defines a $G$-equivariant map $\bigoplus_{x\in X^0}{\mathcal F}(x)\to V$. It factors over the quotient $H_0(X,{\mathcal F})$ of $\bigoplus_{x\in X^0}{\mathcal F}(x)$ as follows immediately from the compatibility requirements on morphisms of coefficient systems (and from the alternating sign inserted into the transition maps for ${\mathcal F}_V$). It is clear that we have defined two assignments which are inverse to each other.\hfill$\Box$\\ 

\section{$p$-torsion representations of ${\rm SL}_2({\mathbb F}_p)$}

Let $\overline{G}={\rm SL}_2({\mathbb F}_p)$, let
$\overline{N}\ne\overline{N}'$ be the unipotent radicals of two opposite
(equivalently: two different) Borel
subgroups in $\overline{G}$. For $\Lambda\in{\rm Art}({\mathfrak o})$ let $\overline{\mathcal J}_{\Lambda}={\rm ind}_{\overline{N}}^{\overline{G}}{\bf 1}_{\Lambda}$ denote the $\overline{G}$-representation on the ${\Lambda}$-module of ${\Lambda}$-valued functions on $\overline{N}\backslash \overline{G}$.\\

\begin{lem}\label{herzjesu} Let $W$ be a $k[\overline{G}]$-module which is
  generated by $W^{{\overline{N}'}}$, let $$\eta:k[\overline{N}]\otimes_k
  W^{{\overline{N}'}} \longrightarrow W$$ denote the morphism of
  $k[\overline{N}]$-modules induced from the inclusion $W^{{\overline{N}'}}\to
  W$, let $$\epsilon: k[\overline{N}]\otimes_k W^{\overline{N}'}\to
  W^{\overline{N}'},\quad\quad \bar{n}\otimes w\mapsto w\mbox{ for }\bar{n}\in \overline{N}$$ denote the
  augmentation map. Then:\\(i) $\eta$ is surjective, i.e. $W$ is generated by $W^{{\overline{N}'}}$ even as a $k[\overline{N}]$-module, and\\(ii) ${\rm ker}(\eta)\subset{\rm ker}(\epsilon)$.\\(iii) The map $H^1({\mathbb Z}_p,\eta):H^1({\mathbb Z}_p,k[\overline{N}]\otimes_k W^{{\overline{N}'}})\to H^1({\mathbb Z}_p,W)$ (continuous cohomology) is bijective for every surjection ${\mathbb Z}_p\to\overline{N}$ .
\end{lem}

{\sc Proof:} (i) Suppose that we are given a surjective map $Y\to W$ of $k[\overline{G}]$-modules. Looking at the commutative diagram$$\xymatrix{k[\overline{N}]\otimes_kY^{\overline{N}'}\ar[d]\ar[r]& Y\ar[d]\\ k[\overline{N}]\otimes_kW^{\overline{N}'}\ar[r]&W}$$we see that statement (i) for the  $k[\overline{G}]$-module $Y$ follows from statement (i) for the  $k[\overline{G}]$-module $W$. Also observe that if we know statement (i) for the direct summands of a direct sum of $k[\overline{G}]$-modules, then we obtain statement (i) also for this direct sum.

Now any $k[\overline{G}]$-module generated by its ${{\overline{N}'}}$-invariants is
also generated by its ${{\overline{N}}}$-invariants (since ${{\overline{N}}}$
and ${{\overline{N}'}}$ are conjugate), hence it is a quotient of a direct sum
of copies of $\overline{\mathcal J}_k$. From \cite{vigcomp} (which is
formulated for ${\rm GL}_2({\mathbb F}_p)$
instead of ${\rm SL}_2({\mathbb F}_p)$) we infer: $\overline{\mathcal J}_k$ is
a direct sum of principal series representations (over $k$) of
$\overline{G}$, and each such principal series representation has length 2 and a
$2$-dimensional subspace of ${\overline{N}}$-invariants. 

Together with our initial remarks we see that to prove statement (i) we may assume that $W$ has length $2$ and that ${\rm dim}_k W^{{\overline{N}}}=2$ (by the above argument we may even assume that $W$ is a principal series representation, but we do not see if this additional assumption could be used to shorten the following argument).

Let ${\mathfrak N}$ denote the set of unipotent radicals of Borel subgroups in $\overline{G}$ which are different from $\overline{N}$.

{\it Claim: We have $k[\overline{N}].W^{{\overline{N}'}}=\sum_{{\overline{N}''}\in {\mathfrak N}}W^{{\overline{N}''}}$, and this contains a $k[\overline{G}]$-submodule $Z\ne0$.}

For any ${\overline{N}''}\in{\mathfrak N}$ we find some $n''\in \overline{N}$
with $n'' {\overline{N}'}(n'')^{-1}={\overline{N}''}$ and hence with
$n''W^{{\overline{N}'}}=W^{{\overline{N}''}}$; conversely, each $n''
{\overline{N}'}(n'')^{-1}$ with $n''\in \overline{N}$ belongs to ${\mathfrak N}$. This shows $k[\overline{N}].W^{{\overline{N}'}}=\sum_{{\overline{N}''}\in {\mathfrak N}}W^{{\overline{N}''}}$. 

Next notice that
$k[\overline{N}].W^{{\overline{N}'}}=\sum_{{\overline{N}''}\in {\mathfrak N}}W^{{\overline{N}''}}$
is $\overline{N}$-stable, hence must contain a non-zero vector $v$ fixed under
$\overline{N}$ (since $\overline{N}$ is a pro-$p$-group and ${\rm
  char}(k)=p$). Let $\overline{P}$ denote the Borel subgroup containing
${\overline{N}}$. Then $\sum_{{\overline{N}''}\in{\mathfrak N}}W^{{\overline{N}''}}$ is stable
  under $\overline{P}$, since the conjugation action of $\overline{P}$ fixes
  the set ${\mathfrak N}$. Therefore $k[\overline{P}]v$ is contained
  in $\sum_{{\overline{N}''}\in {\mathfrak N}}W^{{\overline{N}''}}$. On the other hand, for $\overline{g}\in\overline{G}$ with
$\overline{g}k[\overline{P}].v\ne k[\overline{P}].v$ we have
$\overline{g}{\overline{N}}\overline{g}^{-1}\in{\mathfrak N}$ (look e.g. at
the Bruhat decomposition). But
$k[\overline{P}].v\subset W^{{\overline{N}}}$ implies
$\overline{g}k[\overline{P}].v\subset
W^{\overline{g}{\overline{N}}\overline{g}^{-1}}$. It follows that $Z=\sum_{\overline{g}\in\overline{G}}\overline{g}k[\overline{P}].v$ is a $\overline{G}$-stable submodule of $\sum_{{\overline{N}''}\in {\mathfrak N}}W^{{\overline{N}''}}$. The claim is proved.

{\it Claim: $W=\sum_{{\overline{N}''}\in {\mathfrak N}}W^{{\overline{N}''}}$}

Since $W$ has length $2$ we are done if also $Z$ has
length $2$. Otherwise $W/Z$ and $Z$ are non-zero irreducible, hence
$W/Z=\sum_{{\overline{N}''}\in {\mathfrak N}}(W/Z)^{{\overline{N}''}}$ and
$Z=\sum_{{\overline{N}''}\in {\mathfrak N}}Z^{{\overline{N}''}}$ (the
respective right hand
sides are $\overline{N}$-stable, hence must contain a non zero $\overline{N}$-stable
vector; by $\overline{G}$-irreducibility this must be a generator of the space of ${{\overline{N}}}$
invariants). But all the maps $W^{{\overline{N}''}}\to
(W/Z)^{{\overline{N}''}}$ are surjective (as the $k[\overline{G}]$-module $W$
has length $2={\rm dim}_k W^{{\overline{N}}}={\rm dim}_k W^{{\overline{N}''}}$). We deduce that
$W=\sum_{{\overline{N}''}\in {\mathfrak N}}W^{{\overline{N}''}}$ as desired.

Since as note above, for any ${\overline{N}''}\in{\mathfrak N}$ we find some $n''\in \overline{N}$
with $n''W^{{\overline{N}'}}=W^{{\overline{N}''}}$ it is clear that $\sum_{{\overline{N}''}\in {\mathfrak N}}W^{{\overline{N}''}}$ is the $k[\overline{N}]$-submodule of $W$ generated by $W^{{\overline{N}'}}$. Statement (i) is proven.     

(ii) {\it Claim: As a $k[\overline{N}]$-module, $W$ can not be generated by fewer than $\dim_k(W^{{\overline{N}'}})$ many elements.} 

 Any (finitely generated) $k[\overline{N}]$-module admits a direct sum decomposition with summands isomorphic to quotients of $k[\overline{N}]$: this can be seen  e.g. by applying the structure theorem for modules over the polynomial ring in one variable over $k$, of which $k[\overline{N}]$ is a quotient. From this we see that the minimal number of elements needed to generate such a $k[\overline{N}]$-module is the same as the dimension of its space of $\overline{N}$-invariants. As ${\rm dim}_k(W^{{\overline{N}}})=\dim_k(W^{{\overline{N}'}})$ the claim follows.

Now let $x\in\ker(\eta)$. Choose a finitely generated subspace $S$ of $W^{{\overline{N}'}}$ with $x\in k[\overline{N}]\otimes_kS$, let $W'$ be the $k[\overline{G}]$-module generated by $\eta(k[\overline{N}]\otimes_kS)$: replacing $W$ by $W'$ we see that we may assume that $W$ is finitely generated. If $x\notin \ker(\epsilon)$ then the class of $x$ in $(k[\overline{N}]\otimes_k W^{{\overline{N}'}})\otimes_{k[\overline{N}]}k$ does not vanish (here $k[\overline{N}]\to k$ is the augmentation, its kernel is the maximal ideal of the local ring $k[\overline{N}]$). Therefore, by Nakayama's Lemma, the quotient $(k[\overline{N}]\otimes_k W^{{\overline{N}'}})/(k[\overline{N}].x)$ can be generated by fewer than $\dim_{k}(W^{{\overline{N}'}})$ many elements. As $x\in \ker(\eta)$ this is a contradiction to the claim stated and proved above.

(iii) Let $\gamma:{\mathbb Z}_p\to k[\overline{N}]\otimes_k
W^{{\overline{N}'}}$ be a $1$-cocycle such that $\eta\circ\gamma:{\mathbb
  Z}_p\to W$ is a coboundary. As $\eta$ is surjective we may modify $\gamma$
by a coboundary such that now $\eta\circ\gamma=0$. Let $c\in {\mathbb Z}_p$ be
a generator. Since $\eta(\gamma(c))=0$ we have $\gamma(c)\in{\rm
  ker}(\epsilon)$ by (ii). But ${\rm
  ker}(\epsilon)=(c-1)k[\overline{N}]\otimes_k W^{{\overline{N}'}}$, so
$\gamma(c)=cf-f$ for some $f\in k[\overline{N}]\otimes_k
W^{{\overline{N}'}}$. Since $c$ generates ${\mathbb Z}_p$ the cocycle
condition on $\gamma$ shows $\gamma(c')=c'f-f$ for any $c'\in {\mathbb Z}_p$,
so $\gamma$ is a coboundary. We have shown injectivity of $H^1({\mathbb
  Z}_p,\eta)$. The surjectivity of $H^1({\mathbb Z}_p,\eta)$ follows from the
surjectivity of $\eta$ and from $H^2({\mathbb Z}_p,?)=0$ (for the latter fact
see e.g. \cite{sega} chapter I, section 3.4, Corollary).\hfill$\Box$\\ 

{\it Remark:} In the setting of Lemma \ref{herzjesu}(iii) let $\overline{c}\in
\overline{N}$ be the image of a chosen generator $c$ of ${\mathbb
  Z}_p$. Shapiro's Lemma provides an isomorphism$$W^{{\overline{N}'}}\cong{\rm
  Hom}(p{\mathbb Z}_p,W^{{\overline{N}'}})\cong H^1({\mathbb
  Z}_p,k[\overline{N}]\otimes_k W^{{\overline{N}'}})$$sending $x\in
W^{{\overline{N}'}}$ to the class of the cocycle ${\mathbb Z}_p\to
k[\overline{N}]\otimes_k W^{{\overline{N}'}}$ which maps $c$ to
$\overline{c}\otimes x$. Composing with the isomorphism $H^1({\mathbb
  Z}_p,k[\overline{N}]\otimes_k W^{{\overline{N}'}})\cong H^1({\mathbb
  Z}_p,W)$ of Lemma \ref{herzjesu} (iii)
and with the natural isomorphism $H^1({\mathbb
  Z}_p,W)\cong W_{{\overline{N}}}$ we obtain the composite of natural maps
$W^{{\overline{N}}'}\to W\to W_{{\overline{N}}}$, which therefore is an isomorphism. In particular, on the category of $k[\overline{G}]$-modules generated by their ${\overline{N}}'$-invariants the functors of taking ${\overline{N}}'$-invariants resp. ${\overline{N}}$-coinvariants are equivalent, and in particular are both exact.

\begin{lem}\label{qpfpspec} Let $\Lambda\in{\rm Art}({\mathfrak o})$. (i) Let $V\to W$ be a surjective map of $\Lambda[\overline{G}]$-modules which are generated by their $\overline{N}$-invariant vectors. Then the induced map $V^{\overline{N}}\to W^{\overline{N}}$ is surjective, too.\\(ii) Let $V\to W$ be an injective map of $\Lambda[\overline{G}]$-modules and suppose that $W$ is generated by its $\overline{N}$-invariant vectors. Then the same is true for $V$.
\end{lem}

{\sc Proof:} (i) We may lift any element in $W^{\overline{N}}$ to $V$,
express it as a linear combination of ${\overline{N}}'$-invariant elements and
consider the $\Lambda[\overline{G}]$-sub module inside $V$ generated by these
${\overline{N}}'$-invariant element --- we may
therefore assume that $V$ is finitely generated. By Nakayama's Lemma it is
enough to prove that $V^{\overline{N}}\to W^{\overline{N}}\otimes_{\Lambda} k$ is
surjective. Therefore we may assume that $W$ is in fact a
$k[\overline{G}]$-module. We argue by induction on the minimal number of
generators for $V$. Choose a system $v_1,\ldots,v_s$ in $V^{\overline{N}}$
which generates $V$ as a ${\Lambda}[\overline{G}]$-module and with minimal
$s$. Let $V'\subset V$ be the ${\Lambda}[\overline{G}]$-submodule generated by
$v_1,\ldots,v_{s-1}$, let $W'\subset W$ be its image. We have a surjective
morphism of $\Lambda[\overline{G}]$-modules $V'\oplus \overline{\mathcal
  J}_{\Lambda}\to V$. By induction hypothesis, $(V')^{\overline{N}}\to
(W')^{\overline{N}}$ is surjective. Therefore it is enough to show that the
chain of surjections of $\Lambda[\overline{G}]$-modules $$\overline{\mathcal
  J}_{\Lambda}\longrightarrow \overline{\mathcal J}_{k}\longrightarrow
W/W'$$remains surjective on $\overline{N}$-invariants. For the first arrow
this is clear, it remains to check this for the second arrow. In other words,
we need to check the Lemma for $\Lambda=k$ and $V=\overline{\mathcal
  J}_{k}$. As such it becomes a special case of Lemma \ref{herzjesu}.
(resp. its proof, resp. the remark following it: since for any irreducible
$k[\overline{G}]$-module the space of $\overline{N}$-invariants is
one-dimensional, it is enough to observe that the length of $\overline{\mathcal
  J}_{k}$ equals $\dim_k(\overline{\mathcal J}_{k}^{\overline{N}})$: both
  numbers are equal to $2$).\\(ii) Clearly we may assume that $W$ is a finitely
  generated $\Lambda[\overline{G}]$-module.

{\it Step 1:} First consider the case where $W$
  is actually a $k[\overline{G}]$-module. Let $V'\subset V$ be the
  $k[\overline{G}]$-module generated by $V^{\overline{N}}$. By (i) the
  projection $W\to W/V'$ remains surjective on $\overline{N}$-invariants, in
  particular $V^{\overline{N}}\to (V/V')^{\overline{N}}$ is surjective.  But
  $(V/V')^{\overline{N}}\ne 0$ if $V/V'\ne 0$ since $\overline{N}$ is a
  pro-$p$-group and $V/V'$ is an ${\mathbb F}_p$-vector space. We see $V=V'$.

{\it Step 2:} In the general case we argue by induction on the ${\mathfrak
  o}$-length of $\Lambda$. If $\Lambda\to k$ is not bijective choose a
non-zero element $\pi$ in its kernel with $\pi^2=0$. Then
$\overline{\Lambda}={\Lambda}/(\pi)\in{\rm Art}({\mathfrak o})$ has smaller
${\mathfrak o}$-length than $\Lambda$. By induction hypothesis, the image of $V$ in $W\otimes \overline{\Lambda}$ is
generated by its $\overline{N}$-invariants. By statement (i) we may lift these
$\overline{N}$-invariants to $\overline{N}$-invariants in $V$. Letting
$V'={\rm ker}[V\to W\otimes \overline{\Lambda}]$ we therefore have $V'+\Lambda[\overline{G}]V^{\overline{N}}=V$. It suffices to
prove that $V'$ is generated by its $\overline{N}$-invariants. This follows again
from the induction hypothesis, now applied to the injection $V'\to{\rm ker}[W\to W\otimes \overline{\Lambda}]$,
since $\pi W={\rm ker}[W\to W\otimes\overline{\Lambda}]$ is in fact a
$\overline{\Lambda}[\overline{G}]$-module.\hfill$\Box$\\

\section{Representations versus coefficient systems}

For $x\in X_+^0$ let $$\Theta(x)=\{\tau\in X^1_+\,|\,x\in\tau\mbox{ and
}\tau\mbox{ is not fixed by } N_0\cap K_x\}.$$(Notice that $N_0\cap K_x=N\cap K_x$.) Thus $\Theta(x)$ is the set of all
$1$-simplices $\tau$ containing $x$ except of the single one whose second vertex lies on the
geodesic from $x_-$ to $x$. 

Let ${\mathcal F}$ be an $N_0$-equivariant coefficient system of abelian groups on $X_+$ with injective transition maps (i.e. ${\mathcal F}(\tau)\to {\mathcal F}(x)$ is injective for $x\in\tau\in X_+^1$).

\begin{lem}\label{cogtri} Suppose that for each $x\in X_+^0$ the map
  \begin{gather}H^1(N_0\cap K_x,\bigoplus_{\tau\in\Theta(x)}{\mathcal F}(\tau))\to H^1(N_0\cap K_x,{\mathcal
    F}(x))\label{loccrit}\end{gather} induced by the natural map $\oplus_{\tau\in\Theta(x)}{\mathcal
    F}(\tau)\to{\mathcal F}(x)$ is injective. Then the natural map $C_0(X_+,{\mathcal F})^{N_0}\to H_0(X_+,{\mathcal F})^{N_0}$ is surjective.
\end{lem} 

{\sc Proof:} We write $C_j(X_+,{\mathcal F})=\oplus_{\tau\in
  X^j_+}{\mathcal F}(\tau)$ for the group of $j$-chains on $X_+$, for $j\in\{0,1\}$. As
${\mathcal F}$ has injective transition maps the complex $$0\longrightarrow
C_1(X_+,{\mathcal F})\longrightarrow C_0(X_+,{\mathcal F})\longrightarrow H_0(X_+,{\mathcal F})\longrightarrow0$$
is exact. In cohomology we obtain the exact sequence$$C_0(X_+,{\mathcal
  F})^{N_0}\longrightarrow H_0(X_+,{\mathcal F})^{N_0}\longrightarrow
H^1(N_0,C_1(X_+,{\mathcal F}))\longrightarrow H^1(N_0,C_0(X_+,{\mathcal
  F}))$$so we need to show that $H^1(N_0,C_1(X_+,{\mathcal F}))\to
H^1(N_0,C_0(X_+,{\mathcal F}))$ is injective. For $j\in\{0,1\}$ let $A^j_+$ denote the set of $\tau\in X^j_+$ contained in the apartment of $X$ corresponding to
$T$. Then $X^j_+=\coprod_{\tau\in A^j_+}N_0.\tau$ (disjoint union) and hence$$C_j(X_+,{\mathcal F})=\bigoplus_{\eta\in A^j_+} \bigoplus_{\tau\in
  N_0.\eta}{\mathcal F}(\tau)=\bigoplus_{\eta\in A^j_+}{\rm ind}_{N_0(\eta)}^{N_0}{\mathcal F}(\eta)$$for $j\in\{0,1\}$,
where $N_0(\eta)$ denotes the stabilizer of $\eta$ in $N_0$. Thus,$$C_0(X_+,{\mathcal F})=\bigoplus_{x\in A^0_+} \bigoplus_{y\in
  N_0.x}{\mathcal F}(y)=\bigoplus_{x\in A^0_+}{\rm ind}_{N_0\cap K_x}^{N_0}{\mathcal
  F}(x),$$$$C_1(X_+,{\mathcal F})=\bigoplus_{x\in A^0_+}\bigoplus_{y\in
  N_0.x}\bigoplus_{\tau\in \Theta(y)}{\mathcal
  F}(\tau)=\bigoplus_{x\in A^0_+}{\rm ind}_{N_0\cap K_x}^{N_0}(\bigoplus_{\tau\in\Theta(x)}{\mathcal F}(\tau)).$$By
Shapiro's Lemma we obtain$$H^1(N_0,C_0(X_+,{\mathcal F}))=\bigoplus_{x\in
  A^0_+}H^1(N_0\cap K_x,{\mathcal F}(x)),$$$$H^1(N_0,C_1(X_+,{\mathcal F}))=\bigoplus_{x\in
  A^0_+}H^1(N_0\cap K_x,\bigoplus_{\tau\in\Theta(x)}{\mathcal
  F}(\tau)).$$Therefore we conclude with the injectivity of the
maps (\ref{loccrit}) for each $x\in A^0_+$. \hfill$\Box$\\

\begin{satz}\label{corrpro} For any ${\mathcal F}\in{\mathfrak C}(k)$ the natural map ${\mathcal F}(\sigma)\to H_0({X}_+,{\mathcal F})^{N_{0}}$ is bijective.
\end{satz}

{\sc Proof:} The map in question naturally factors as\begin{gather}{\mathcal F}(\sigma)\longrightarrow\frac{C_0({X}_+,{\mathcal F})^{N_{0}}}{C_1({X}_+,{\mathcal
  F})^{N_{0}}}\longrightarrow H_0({X}_+,{\mathcal
F})^{N_{0}}.\label{factor}\end{gather}The first arrow in (\ref{factor}) is bijective by condition (b) in the definition of ${\mathfrak
  C}(k)$. Now let $x\in X_+^0$. Then $\overline{N}=N_0\cap K_x/(N_0\cap K_x)^p$ is the unipotent radical of a
Borel subgroup in the reductive quotient $\overline{G}$ of $K_x$. Let
$\overline{N}'$ denote the unipotent radical of another Borel subgroup in
$\overline{G}$; it fixes an element of $\Theta(x)$, while $\overline{N}$ acts
simply transitively on $\Theta(x)$. Putting $W={\mathcal
  F}(x)$ we therefore obtain an $\overline{N}$-equivariant identification $k[\overline{N}]\otimes
W^{\overline{N}'}=\bigoplus_{\tau\in\Theta(x)}{\mathcal F}(\tau)$, and condition (c) in the definition of ${\mathfrak
  C}(k)$ shows that we are in the
setting of Lemma \ref{herzjesu}. We have $N_0\cap K_x\cong {\mathbb Z}_p$, so Lemma
\ref{herzjesu} tells us that the
map (\ref{loccrit}) is injective. Therefore the second arrow in (\ref{factor}) is bijective by Lemma \ref{cogtri}. \hfill$\Box$\\ 

\begin{satz}\label{h0} Let $\Lambda\in{\rm Art}({\mathfrak o})$.\\(i) For $V\in{\mathfrak R}(\Lambda)$ the natural map $H_0(X,{\mathcal F}_V)\to V$ is an isomorphism.\\(ii) For ${\mathcal F}\in{\mathfrak C}(\Lambda)$ the natural map ${\mathcal F}\to{\mathcal F}_{H_0({X},{\mathcal F})}$ is an isomorphism.\\(iii) ${\mathfrak R}(\Lambda)$ and ${\mathfrak C}(\Lambda)$ are abelian categories. For any epimorphism $\Lambda'\to\Lambda$ in ${\rm Art}({\mathfrak o})$, any $W\in {\mathfrak R}(\Lambda')$ and $V\in{\mathfrak R}(\Lambda)$, if $W\to V$ is a surjective map in ${\mathfrak R}(\Lambda')$ then the induced map $W^{U_{\sigma}}\to V^{U_{\sigma}}$ is surjective.
\end{satz}

{\sc Proof:} {\it Step 1: The proof of statement (i) for $V\in{\mathfrak
    R}(k)$}

Consider the composite $$V^{U_{\sigma}}={\mathcal F}_V({\sigma})\longrightarrow
H_0({X}_+,{\mathcal F}_V)^{N_{0}}\longrightarrow V.$$ The first arrow is
bijective by Theorem \ref{corrpro}, therefore the second one is injective, and
as $k$ has characteristic $p$ and ${N_{0}}$ is a pro-$p$-group it follows
that $H_0({X}_+,{\mathcal F}_V)\to V$ is injective. Now $X$ can be written as
an increasing union $X={\cup}_iX_{+,i}$ of half trees $X_{+,{i}}$ contained
in $X$, all of them isomorphic with $X_+$. On each of them we may repeat the
same constructions as for $X_+$. It follows that $H_0({X},{\mathcal
  F}_V)=\cup_i H_0({X}_{+,i},{\mathcal F}_V)\to V$ is injective, and hence bijective.

{\it Step 2: Statement (i) for a fixed $\Lambda$ implies statement (ii) for
  the same $\Lambda$.} 

To see this let ${\mathcal F}\in{\mathfrak C}(\Lambda)$. As ${\mathcal F}$ has
injective transition maps the natural map ${\mathcal F}(\tau)\to
H_0(X,{\mathcal F})$ is injective for any simplex $\tau$. We may therefore
regard ${\mathcal F}(\tau)$ and ${\mathcal F}_{H_0(X,{\mathcal F})}(\tau)$ as
being contained in $H_0(X,{\mathcal F})$, so our hypotheses on ${\mathcal F}$
and the definition of ${\mathcal F}_{H_0(X,{\mathcal F})}$ show that we may
regard ${\mathcal F}$ as a subcoefficient system of ${\mathcal
  F}_{H_0(X,{\mathcal F})}$. Namely, as $U_{\tau}$ for $\tau\in X^1$ acts trivially on
${\mathcal F}(\tau)$ we have the inclusion maps
$\alpha_{\tau}:{\mathcal F}(\tau)\to H_0(X,{\mathcal F})^{U_{\tau}}={\mathcal
  F}_{H_0(X,{\mathcal F})}(\tau)$, and as ${\mathcal F}(x)=\sum_{\tau\in
  X^1\atop x\in \tau}{\mathcal F}(\tau)$ the $\alpha_{\tau}$ also induce maps $\alpha_{x}:{\mathcal F}(x)=\sum_{\tau\in
  X^1\atop x\in \tau}{\mathcal F}(\tau)\to \sum_{\tau\in
  X^1\atop x\in \tau}H_0(X,{\mathcal F})^{U_{\tau}}={\mathcal
  F}_{H_0(X,{\mathcal F})}(x)$ for $x\in X^0$. In particular we obtain a map
$\alpha:H_0(X,{\mathcal F})\to H_0(X,{\mathcal F}_{H_0(X,{\mathcal F})})$. On the other hand, by the definition of ${\mathcal F}_V$ for
$V=H_0(X,{\mathcal F})$ there is a natural map $\beta: H_0(X,{\mathcal
  F}_{H_0(X,{\mathcal F})})\to H_0(X,{\mathcal F})$. We claim that
$\beta\circ\alpha$ is the identity on $H_0(X,{\mathcal F})$. Indeed, as
$H_0(X,{\mathcal F})$ is generated by the images of the natural maps
$\iota_x:{\mathcal F}(x)\to H_0(X,{\mathcal F})$ for all $x\in X^0$, it
is enough to show $\beta\circ\alpha\circ\iota_x=\iota_x$ for all $x\in
X^0$. If $\eta_x:{\mathcal
  F}_{H_0(X,{\mathcal F})}(x)\to H_0(X,{\mathcal
  F}_{H_0(X,{\mathcal F})})$ denotes the natural map, then we have
  $\alpha\circ\iota_x=\eta_x\circ \alpha_x$ by the definition of $\alpha$. Now
  by the definition of $\beta$ we have that $\beta\circ \eta_x$ is just the inclusion ${\mathcal
  F}_{H_0(X,{\mathcal F})}(x)=\sum_{\tau\in
  X^1\atop x\in \tau}H_0(X,{\mathcal F})^{U_{\tau}}\to H_0(X,{\mathcal F})$
for all $x\in X^0$. It follows that $\beta\circ\alpha\circ\iota_x=\beta\circ
\eta_x\circ \alpha_x$ is the inclusion of ${\mathcal F}(x)$ into $H_0(X,{\mathcal F})$, i.e. the
map $\iota_x$, as desired. The claim is proven. 

By statement (i), applied to $V=H_0(X,{\mathcal F})$, the map $\beta$ is an isomorphism, hence so is $\alpha$. In particular, $H_0(X,{\mathcal F}_{H_0({X},{\mathcal F})}/{\mathcal F})=0$. But it follows from our hypotheses on ${\mathcal F}$ that for all $x\in\eta\in X^1$ we have $${\mathcal F}_{H_0({X},{\mathcal F})}(\eta)\cap{\mathcal F}(x)={\mathcal F}(\eta)$$inside ${\mathcal F}_{H_0({X},{\mathcal F})}(x)$, i.e. $({\mathcal F}_{H_0({X},{\mathcal F})}/{\mathcal F})(\eta)\to({\mathcal F}_{H_0({X},{\mathcal F})}/{\mathcal F})(x)$ is injective, i.e. the quotient system ${\mathcal F}_{H_0({X},{\mathcal F})}/{\mathcal F}$ has injective transition maps. These two facts together imply ${\mathcal F}_{H_0({X},{\mathcal F})}/{\mathcal F}=0$, i.e. ${\mathcal F}_{H_0({X},{\mathcal F})}={\mathcal F}$.

{\it Step 3: Statements (i) and (ii) for a fixed $\Lambda$ imply statement (iii)
  for the same $\Lambda$.}

Statement (i) and (ii) (and Lemma \ref{adjoin}) together say that the functors $V\mapsto {\mathcal F}_V$ and ${\mathcal F}\mapsto H_0(X,{\mathcal F})$ set up an equivalence between the categories ${\mathfrak R}(\Lambda)$ and ${\mathfrak C}(\Lambda)$. That ${\mathfrak C}(\Lambda)$ is an abelian category follows from Lemma \ref{qpfpspec}. As ${\mathcal F}\mapsto H_0(X,{\mathcal F})$ is exact on ${\mathfrak C}(\Lambda)$ it follows that also ${\mathfrak R}(\Lambda)$ is an abelian category and that $V\mapsto {\mathcal F}_V$ is exact on ${\mathfrak R}(\Lambda)$. In particular, for a surjection $W\to V$ in ${\mathfrak R}(\Lambda)$ the induced map $W^{U_{\sigma}}\to V^{U_{\sigma}}$ is surjective. Let now more generally $W\to V$ be a surjection in ${\mathfrak R}(\Lambda')$ for some epimorphism $\Lambda'\to\Lambda$, with $V\in {\mathfrak R}(\Lambda)$. As $W$ is generated by $W^{U_{\sigma}}$ there is a free $\Lambda'$-module $T$ and, if we endow $T$ with the trivial $U_{\sigma}$-action, a surjection ${\rm ind}_{U_{\sigma}}^GT\to W$. Its composite with $W\to V$ can alternatively be factored as $${\rm ind}_{U_{\sigma}}^GT\longrightarrow{\rm ind}_{U_{\sigma}}^G(T\otimes_{\Lambda'}\Lambda)\longrightarrow V.$$The first one of these arrows remains surjective on $U_{\sigma}$-invariants for general reasons, for the second one this is true by what we observed half a minute ago.

{\it Step 4: Statement (i) for a given $\Lambda\in {\rm
    Art}({\mathfrak o})$ follows from statements (i), (ii) and (iii) for all elements
  in ${\rm Art}({\mathfrak o})$ of smaller ${\mathfrak o}$-length.}

 If $\Lambda\to k$ is not bijective choose a non-zero element $\pi$ in its
 kernel with $\pi^2=0$. Then $\overline{\Lambda}={\Lambda}/(\pi)\in{\rm
   Art}({\mathfrak o})$ has smaller ${\mathfrak o}$-length than
 $\Lambda$. For $V\in{\mathfrak R}(\Lambda)$ the representations
 $\pi V$ and $V\otimes \overline{\Lambda}$ are in fact
 $\overline{\Lambda}$-modules and hence lie in ${\mathfrak
   R}(\overline{\Lambda})$. From statement (iii), applied to $V\to V\otimes
 \overline{\Lambda}$ and $\pi:V\to \pi V$, it follows that the natural maps
 ${\mathcal F}_V\otimes \overline{\Lambda}\to{\mathcal F}_{V\otimes
   \overline{\Lambda}}$ and $\pi({\mathcal F}_V)\to{\mathcal F}_{\pi V}$ are
 surjective at $1$-simplices, hence also at $0$-simplices. On the other hand,
 for formal reasons these maps are also injective at all simplices, hence are
 isomorphisms of coefficient systems. Therefore, the exact sequence of coefficient
 systems $$0\longrightarrow\pi({\mathcal F}_V)\longrightarrow{\mathcal F}_V\longrightarrow {\mathcal F}_V\otimes\overline{\Lambda}\longrightarrow0$$ can be read as the following sequence of coefficient
 systems:$$0\longrightarrow{\mathcal F}_{\pi V}\longrightarrow{\mathcal F}_V\longrightarrow {\mathcal F}_{V\otimes
 \overline{\Lambda}}\longrightarrow0.$$We obtain the commutative diagram 

$$\xymatrix{0\ar[r]& H_0(X,{\mathcal F}_{\pi V}) \ar[d]\ar[r]& H_0(X,{\mathcal F}_{V})\ar[d]\ar[r]& H_0(X,{\mathcal F}_{V\otimes
 \overline{\Lambda}})\ar[d]\ar[r]&0\\0\ar[r]& \pi V \ar[r]&V\ar[r]&V\otimes
 \overline{\Lambda}
       \ar[r]&0}$$with exact rows, observing $H_1(X,{\mathcal
   F}_{V\otimes\overline{\Lambda}})=0$ (which holds true as ${\mathcal F}\otimes
 \overline{\Lambda}$ has injective transition maps). Applying the induction hypothesis to
 $\pi V$ and $V\otimes \overline{\Lambda}$ we conclude.

 {\it Step 5:} The combination of all these steps gives the Theorem.\hfill$\Box$\\

{\bf Definition:} For $\Lambda\in{\rm Art}({\mathfrak o})$ let ${\mathfrak R}_{\ast}(\Lambda)$ denote the category of smooth ${\rm GL}_2({\mathbb Q}_p)$-representations on $\Lambda$-modules which, when restricted to $G={\rm SL}_2({\mathbb Q}_p)$, belong to ${\mathfrak R}(\Lambda)$. Let ${\mathfrak C}_{\ast}(\Lambda)$ denote the category of ${\rm GL}_2({\mathbb Q}_p)$-equivariant coefficient systems on $X$ which, when restricted to $G={\rm SL}_2({\mathbb Q}_p)$, belong to ${\mathfrak C}(\Lambda)$. 

\begin{kor}\label{main} Let $\Lambda\in{\rm Art}({\mathfrak o})$, let $?$ be the empty symbol or $?={\ast}$. The categories ${\mathfrak R}_?(\Lambda)$ and ${\mathfrak C}_?(\Lambda)$ are abelian. The functors $${\mathfrak R}_?(\Lambda)\longrightarrow{\mathfrak C}_?(\Lambda),\quad\quad V\mapsto {\mathcal F}_V$$ and $${\mathfrak C}_?(\Lambda)\longrightarrow{\mathfrak R}_?(\Lambda),\quad\quad {\mathcal F}\mapsto H_0(X,{\mathcal F})$$ are exact and quasi-inverse to each other.
\end{kor}\hfill$\Box$

\section{Complements}

\subsection{A cohomology-free proof of Theorem \ref{corrpro}.}
\label{secondproof}

Let $d:X^0\times X^0\to{\mathbb Z}_{\ge0}$ denote the
counting-vertices-on-geodesics distance. Let $Y_{-1}=\{x_-\}$. For $m\in{\mathbb Z}_{\ge0}$ let
$Y_m=\{x\in{X}_{+}^0\,|\,d(x,x_+)=m\}$ and $Y_{m,m+1}=\{\{x_1,x_2\}\in
X^1\,|\,x_1\in Y_m, x_2\in Y_{m+1}\}$. Thus, $X_+^0=\coprod_{m\ge0}Y_m$ and $X_+^1=\coprod_{m\ge0}Y_{m,m+1}$. For $m\in{\mathbb Z}_{\ge1}$ let
$N_0^{(m)}$ denote the pointwise stabilizer of $Y_{m-1}$ in $N_0$. Thus
$N_0=N_0^{(1)}$. For $v\in X^0$ let $U_v$ denote the kernel of the surjection
of $K_x$ onto its reductive quotient (which is isomorphic with
${\overline{G}}$). Write $${\mathcal F}(m)=\bigoplus_{z\in Y_m}{\mathcal
   F}(z),\quad\quad\quad{\mathcal F}(m,m+1)=\bigoplus_{\eta\in
   Y_{m,m+1}}{\mathcal F}(\eta).$$

{\sc Second Proof of Theorem \ref{corrpro}:} The injectivity of ${\mathcal F}(\sigma)\to H_0({X}_+,{\mathcal F})^{N_{0}^{(1)}}$ follows from that of ${\mathcal F}(\sigma)\to{\mathcal F}(x_+)$. 

Surjectivity. Let $\partial:C_1(X,{\mathcal F}){\to}C_0(X,{\mathcal F})$ denote the differential from $1$-chains to $0$-chains. For a $0$-chain $c=(c_v)_{v\in X_+^0}$ we write $c(m)=\sum_{v\in Y_m}c_v\in{\mathcal F}(m)$, for a $1$-chain $b=(b_{\tau})_{\tau\in {X}_+^1}$ we write $b(m,m+1)=\sum_{\tau\in Y_{m,m+1}}b_{\tau}\in {\mathcal F}(m,m+1)$. Let $g\in N_0^{(1)}\cong {\mathbb Z}_p$ be a topological generator. Then $N_0^{(m+1)}$ is topologically generated by $g^{p^m}$ for any $m\ge0$, and hence so is its quotient $U_{\{v_1,v_2\}}/{U_{v_1}}$ for any $\{v_1,v_2\}\in X^1$, $v_1\in Y_m$, $v_2\in Y_{m-1}$.        

Let the $0$-chain $c$ represent a non zero class $[c]$ in $H_0({X}_+,{\mathcal F})^{N_{0}^{(1)}}$. Let $n(c)\in{\mathbb Z}_{\ge0}$ be minimal with ${\rm supp}(c)\subset \cup_{n\le n(c)}Y_n$. By induction on $n(c)$ we now show that $[c]$ lies in the image of ${\mathcal F}(\sigma)$, as required. If $n(c)>0$ condition (c) (in the definition of ${\mathfrak C}(k)$) shows that we may assume $c(m)\in {\mathcal F}(m)^{N_{0}^{(m+1)}}$ for all $m<n(c)$. If we knew that $c(n(c))\in{\mathcal F}(n(c))$ is fixed under $N_{0}{(n(c)+1)}$ --- or equivalently: by $g^{p^{n(c)}}$ --- then condition (b) (in the definition of ${\mathfrak C}(k)$) shows that\\ --- if $n(c)=0$ then we are done, and\\--- if $n(c)>0$ then we find a $c'$ with $n(c')=n(c)-1$ and $[c]=[c']$, so our induction hypothesis applies to $c'$. 

Thus, it remains to show $g^{p^{n(c)}}(c(n(c)))=c(n(c))$. As $[c]$ is $N_{0}^{(1)}$-fixed there is a $1$-chain $b$ such that $gc=c+\partial(b)$, hence $$g^sc=c+\sum_{i=0}^{s-1}\partial(g^ib)\quad\quad\mbox{ for any }s\in{\mathbb N}.$$By induction on $m$ we show that $\sum_{i=0}^{p^{m+1}-1}g^i{b}(m,m+1)=0$ for any $0\le m<n(c)$. Fix $0\le m<n(c)$, and if even $0< m<n(c)$ then assume $\sum_{i=0}^{p^{m}-1}g^ib(m-1,m)=0$. Since we have $g^{p^{m}}c(m)=c(m)$ by assumption this says$$0=(g^{p^{m}}c-c)(m)=(\sum_{i=0}^{p^{m}-1}\partial(g^i{b}))(m)=\gamma_m^{m+1}(\sum_{i=0}^{p^{m}-1}g^i{b})(m,m+1)$$where $\gamma_m^{m+1}:{\mathcal F}(m,m+1)\to{\mathcal F}(m)$ is the natural map. Choose a system $\{\eta\}$ in $Y_{m,m+1}$ such that each $v\in Y_m$ is contained in exactly one $\eta$. For each such $v$ and $\eta$ we make the identifications $\overline{G}=K_{v}/U_{v}$ and $\overline{N}=N_{0}^{(m+1)}/N_0^{(m+2)}$, and we let $\overline{N}'\subset\overline{G}=K_{v}/U_{v}$ be the Borel subgroup fixing $\eta$. Together with the natural identification$$k[N_{0}^{(m+1)}/N_0^{(m+2)}]\otimes_k(\bigoplus_{\eta}{\mathcal F}(\eta))\cong{\mathcal F}(m,m+1)$$the map $\gamma_m^{m+1}$ then becomes the induced map as considered in Lemma \ref{herzjesu}, and the summation map $\sum_{j=0}^{p-1}g^{jp^m}(.)$ becomes the augmentation map in Lemma \ref{herzjesu}. Therefore Lemma \ref{herzjesu} tells us that$$(\sum_{i=0}^{p^{m+1}-1}g^i{b})(m,m+1)=\sum_{j=0}^{p-1}g^{jp^m}(\sum_{i=0}^{p^{m}-1}g^i{b})(m,m+1)=0$$and the induction proof is complete. The result for $m=n(c)-1$ says $$(\sum_{i=0}^{p^{n(c)}-1}g^i{b})(n(c)-1,n(c))=0.$$ On the other hand, from (i) we see $b(n(c),n(c)+1)=0$. Together$$g^{p^{n(c)}}(c(n(c)))=c(n(c))+(\partial(\sum_{i=0}^{p^{n(c)}-1}g^i{b}))(n(c))=c(n(c))$$and we are done.\hfill$\Box$

\subsection{Hecke modules}

 Let $\Lambda\in{\rm Art}({\mathfrak o})$. Let ${\mathcal J}_{\Lambda,{\ast}}={\rm ind}_{U_{\sigma}}^{{\rm GL}_2({\mathbb Q}_p)}{\bf 1}_{\Lambda}$ (compact induction). Let ${\mathcal H}_{\Lambda}={\rm End}_{\Lambda[{\rm GL}_2({\mathbb Q}_p)]}({\mathcal J}_{\Lambda,{\ast}})$ be the pro-$p$-Iwahori-Hecke algebra. Let ${\mathfrak M}_{\ast}(\Lambda)$ denote the category of right-${\mathcal H}_{\Lambda}$-modules.

 Let $K_{x_+,{\ast}}$ denote the maximal compact open subgroup in ${\rm GL}_2({\mathbb Q}_p)$ fixing the vertex $x_+$. Its reductive quotient is isomorphic with ${\rm GL}_2({\mathbb F}_p)$, and the subgroup $U_{\sigma}$ of $K_{x_+,{\ast}}$ maps to the unipotent radical of a Borel subgroup $\overline{N}$ in ${\rm GL}_2({\mathbb F}_p)$. Let $$\overline{{\mathcal J}}_{\Lambda,{\ast}}={\rm ind}_{U_{\sigma}}^{K_{x_+,{\ast}}}{\bf 1}_{\Lambda}={\rm ind}_{\overline{N}}^{{\rm GL}_2({\mathbb F}_p)}{\bf 1}_{\Lambda}.$$ We may view the finite Hecke algebra ${\mathcal H}^{\natural}_{\Lambda}={\rm End}_{\Lambda[K_{x_+,{\ast}}]}(\overline{{\mathcal J}}_{\Lambda,\ast})$ as a subalgebra of ${\mathcal H}_{\Lambda}$.

\begin{pro}\label{finflat} Let $\Lambda\in{\rm Art}({\mathfrak o})$. Suppose that $\overline{{\mathcal J}}_{\Lambda,{\ast}}$ is flat as an ${\mathcal H}^{\natural}_{\Lambda}$(-left) module. Then the functor \begin{gather}{\mathfrak R}_{\ast}(\Lambda)\longrightarrow {\mathfrak M}_{\ast}(\Lambda),\quad\quad V\mapsto V^{U_{\sigma}}\label{racor}\end{gather}is exact and an equivalence of categories.
\end{pro}

The action of ${\mathcal H}_{\Lambda}$ on $V^{U_{\sigma}}$ results from the isomorphisms$$V^{U_{\sigma}}={\rm Hom}_{U_{\sigma}}({\bf 1}_{\Lambda},V)={\rm Hom}_{\Lambda[{\rm GL}_2({\mathbb Q}_p)]}({\mathcal J}_{\Lambda,{\ast}},V).$$

{\sc Proof:} To an ${\mathcal H}^{\natural}_{\Lambda}$-right module $M$ we
assign the $\Lambda[K_{x_+,{\ast}}]$-module $${\mathcal
  K}(M)=M\otimes_{{\mathcal H}^{\natural}_{\Lambda}}\overline{{\mathcal
    J}}_{\Lambda,\ast}.$$Let $\varphi$ denote the unique element in $\overline{{\mathcal
    J}}_{\Lambda,\ast}$ supported on $U_{\sigma}$ and taking value
$1\in\Lambda$ there. We then have the natural map \begin{gather}M\longrightarrow {\mathcal K}(M)^{U_{\sigma}},\quad
m\mapsto m\otimes\varphi.\label{vytastra}\end{gather}We claim that the map (\ref{vytastra}) is
bijective, for any $M$. For this we easily reduce to the case where $M$ has
finite length. Arguing by an induction on the length $\ell(M)$ of $M$ we
consider an exact sequence $0\to M_1\to M\to M_2\to 0$ with
an ${\mathcal H}^{\natural}_{\Lambda}$-right module $M_1$ with
$\ell(M_1)<\ell(M)$ and with an irreducible ${\mathcal
  H}^{\natural}_{\Lambda}$-right module $M_2$. It gives rise to the commutative
diagram  
\begin{gather}\xymatrix{0\ar[r]& M_1 \ar[d]\ar[r]& M\ar[d]\ar[r]& M_2\ar[d]\ar[r]&0\\0\ar[r]& {\mathcal
  K}(M_1)^{U_{\sigma}} \ar[r]&{\mathcal
  K}(M)^{U_{\sigma}}\ar[r]&{\mathcal
  K}(M_2)^{U_{\sigma}}
       \ar[r]&0}.\label{hecdia}\end{gather}By our flatness assumption on $\overline{{\mathcal
         J}}_{\Lambda,{\ast}}$ the sequence $0\to{\mathcal
  K}(M_1)\to {\mathcal
  K}(M)\to {\mathcal
  K}(M_2)\to 0$ is exact. As ${\mathcal
  K}(M)^{U_{\sigma}}\to{\mathcal
  K}(M_2)^{U_{\sigma}}$ is surjective by Lemma \ref{qpfpspec} it follows that the bottom row in the diagram
(\ref{hecdia}) is exact. By induction hypothesis, the first vertical arrow is
an isomorphism. Therefore it remains to show that the third vertical arrow is
an isomorphism. In other words, we have reduced our claim to the case where
$M$ is an irreducible ${\mathcal H}^{\natural}_{\Lambda}$-right module. But
then $M$ clearly also is an irreducible ${\mathcal H}^{\natural}_{k}$-right
module. For such $M$ our claim is shown in the (easy) proof of \cite{pa}
Corollary 3.3. 

Let $\widetilde{I}_{\sigma}$ denote the stabilizer of $\sigma$ in ${\rm GL}_2({\mathbb
  Q}_p)$. In the terminology of \cite{pa}, a pair $(E_1,E_0)$ consisting of a
$\Lambda[\widetilde{I}_{\sigma}]$-module $E_1$, a
$\Lambda[K_{x_+,{\ast}}]$-module $E_2$ and a $\Lambda[\widetilde{I}_{\sigma}\cap K_{x_+,{\ast}}]$-equivariant map $E_1\to
E_2$ is called a $\Lambda$-diagram. It is shown in \cite{pa} Theorem 5.17 that we have an equivalence of
categories between $\Lambda$-diagrams and $G$-equivariant coefficient systems
of $\Lambda$-modules on $X$. (In fact, in \cite{pa} only $\Lambda=\overline{\mathbb F}_p$ is
considered, but the arguments immediately generalize to more general $\Lambda$.)

In our setting, if $M\in {\mathfrak M}_{\ast}(\Lambda)$ then, in addition to its ${\mathcal
  H}^{\natural}_{\Lambda}$-module structure, $M$ is endowed with an action of $\widetilde{I}_{\sigma}$. This provides the
pair $(M,{\mathcal K}(M))$ with the structure of a $\Lambda$-diagram. Hence we
may associate to these data in a functorial way a ${\rm GL}_2({\mathbb Q}_p)$-equivariant
coefficient system ${\mathcal F}^M$ on $X$ with ${\mathcal F}^M(x_+)={\mathcal
  K}(M)$ and ${\mathcal F}^M(\sigma)=M$. From $M\cong {\mathcal
  K}(M)^{U_{\sigma}}$ it follows that ${\mathcal F}^M$ belongs to ${\mathfrak
  C}_{\ast}(\Lambda)$. Thus, we obtain a functor\begin{gather}{\mathfrak M}_{\ast}(\Lambda)\longrightarrow{\mathfrak
  C}_{\ast}(\Lambda),\quad\quad M\mapsto {\mathcal F}^M.\label{vyra}\end{gather}
Now
we observe:

(i) Morphisms in ${\mathcal F}_1\to {\mathcal F}_2$ in ${\mathfrak
  C}_{\ast}(\Lambda)$ are uniquely determined on $1$-simplices, i.e. by the
collection of morphisms ${\mathcal F}_1(\tau)\to {\mathcal F}_2(\tau)$ for all
$\tau\in X^1$: this follows from property (c) in the definition of ${\mathfrak
  C}(\Lambda)$. Similarly, if the ${\mathcal F}_1(\tau)\to {\mathcal
  F}_2(\tau)$ are isomorphisms for all
$\tau\in X^1$, then ${\mathcal F}_1\to {\mathcal F}_2$ is an
isomorphism.

(ii) For $V\in {\mathfrak
  R}_{\ast}(\Lambda)$ our constructions provide a natural morphism ${\mathcal
  F}^{V^U}\to{\mathcal F}_V$ in ${\mathfrak
  C}_{\ast}(\Lambda)$ which is an isomorphism on $1$-simplices. By (i) it must
be an isomorphism itself. Thus, the composition of functors $${\mathfrak
  R}_{\ast}(\Lambda)\longrightarrow {\mathfrak M}_{\ast}(\Lambda)\longrightarrow{\mathfrak
  C}_{\ast}(\Lambda),\quad V\mapsto {\mathcal
  F}^{V^U}$$is isomorphic with the functor $V\mapsto {\mathcal
  F}_{V}$. 

(iii) The composition of functors $${\mathfrak M}_{\ast}(\Lambda)\longrightarrow{\mathfrak
  C}_{\ast}(\Lambda)\longrightarrow {\mathfrak
  R}_{\ast}(\Lambda)\longrightarrow {\mathfrak M}_{\ast}(\Lambda),\quad
M\mapsto H_0(X,{\mathcal F}^M)^{U_{\sigma}}$$is isomorphic with the identity
functor on ${\mathfrak M}_{\ast}(\Lambda)$. Indeed, by Theorem \ref{h0} we have ${\mathcal
  F}^M(\sigma)\cong {\mathcal F}_{H_0(X,{\mathcal
  F}^M)}(\sigma)$, and this gives $M\cong{\mathcal
  F}^M(\sigma)\cong H_0(X,{\mathcal F}^M)^{U_{\sigma}}$ as
$\Lambda$-modules. The compatibility with the action of $\widetilde{I}_{\sigma}$ is easily
verified. (We remark that we also see directly from (i) that the functor (\ref{vyra}) is fully
faithful.)

From (ii) and (iii) and Theorem \ref{h0} we deduce that the functor $${\mathfrak M}_{\ast}(\Lambda)\to
{\mathfrak R}_{\ast}(\Lambda), \quad\quad M\mapsto H_0(X,{\mathcal
  F}^M)$$ is quasi inverse to the functor
(\ref{racor}). We also see the exactness of these functors.\hfill$\Box$\\ 

{\it Questions:} (a) Is the flatness hypothesis in Proposition \ref{finflat} fulfilled for all $\Lambda$ ? For $\Lambda=k$ it is, as follows from arguments similar to those used in the proofs of Lemmata \ref{herzjesu} and \ref{qpfpspec} (relying on an explicit investigation of $\overline{\mathcal J}_{k,_{\ast}}$). That ${\mathfrak R}_{\ast}(k)\to {\mathfrak M}_{\ast}(k), V\mapsto V^{U_{\sigma}}$ is exact and an equivalence of categories is the main result from \cite{ol}.

(b) One may ask for an ${\rm SL}_2({\mathbb Q}_p)$-analog of Proposition \ref{finflat}.\\

{\it Acknowledgments: I am very grateful to Peter Schneider for several
  discussions on this note. I thank the anonymous referee for helpful comments
  leading to an improved exposition. I thank the Deutsche Forschungs Gemeinschaft (DFG) as this work was done while I was supported by the DFG as a Heisenberg Professor at the Humboldt University, Berlin.} \\

\begin{flushleft} \textsc{Humboldt-Universit\"at zu Berlin\\Institut f\"ur Mathematik\\Rudower Chaussee 25\\12489 Berlin, Germany}\\ \textit{E-mail address}:
gkloenne@math.hu-berlin.de \end{flushleft} \end{document}